# Paving the Way to Smart Micro Energy Internet: Concepts, Design Principles, and Engineering Practices

Shengwei Mei, *Fellow, IEEE,* Rui Li, Xiaodai Xue, Ying Chen, Qiang Lu, *Fellow, IEEE,* Xiaotao Chen, *Carsten D. Ahrens,* Ruomei Li, and Laijun Chen

*Abstract*—The energy internet is one of the most promising future energy infrastructures that could both enhance energy efficiency and improve its operating flexibility. Analogous to the micro-grid, the micro energy internet emphasizes the distribution level and demand side. This paper proposes concepts and design principles of a smart micro energy internet for accommodating micro-grids, distributed poly-generation systems, energy storage facilities, and associated energy distribution infrastructures. Since the dispatch and control system of the smart micro energy internet is responsible for external disturbances, it should be able to approach a satisfactory operating point while supporting multiple criteria, such as safety, economy, and environmental protection. To realize the vision of a smart micro energy internet, an engineering game theory based energy management system with self-approaching-optimum capability is investigated. Based on the proposed concepts, design principles, and energy management system, this paper presents a prototype of China's first conceptual solar-based smart micro energy internet, established in Qinghai University.

## I. Introduction

INCREASING worldwide energy demand has created a global urgency for high-efficiency and environmentally friendly energy utilization systems. It has been validated that the integration of different energy distribution systems could bring additional flexibility and reliability to system operation and enhance energy efficiency [1]. Thus, the coupling of energy infrastructures including district heating network [2], [3], natural gas network [4], [5], the electric vehicle transportation infrastructure [6], [7] and the power network, all have been investigated extensively to improve system efficiency and reliability.

The integrated energy system [1], energy internet [8], multi-source multi-product system [9], multi-energy system [10], as well as the energy micro-grid [11], all provide the possibility for coupling among multi-carrier energy systems, including electricity, natural gas, and heat. The integration of such infrastructures has been viewed as a promising future approach to realize a safe, cost-effective, and environmentally friendly energy system [12], [13].

The energy hub is one of the key components in a multi-carrier energy system. It provides the capabilities of conversion, transmission, and storage among multiple energy carriers [14]. The energy hub also helps to build natural linkages among traditional independent infrastructures, including power distribution network (PDN), gas distribution network (GDN), district heating network (DHN), cold distribution network (CDN), and electrified transportation network (ETN). Combined heat and power units with thermal energy storage [3], compressed air energy storage systems [15], [16], and concentrating solar power facilities [17] are typical energy hubs that offer power and heat co-generation and storage capabilities. These hubs have been investigated at both transmission and distribution-level coupled energy systems [16]–[18].

Since the energy distribution level simplifies the implementation of the coupled energy system functions, research on the future form of such coupled infrastructures both at the distribution level and on the demand side are of great interest [19]. One such coupled system is the micro energy internet, which emphasizes the integration of multiple energy networks at the distribution level. Specifically, the micro energy internet is a distributed form of energy internet composed of distributed renewable energy sources, energy storage facilities, multi-carrier energy sources, multi-carrier loads, and distribution infrastructures, such as PDN, DHN, GDN, and ETN. Based on this, we have the smart micro energy internet, which is formed by capturing multi-criteria self-approaching-optimum capabilities, as illustrated in [20]. In this work, we further exploit the fundamental concepts and design principles of a smart micro energy internet.

Although extra reliability and flexibility can be achieved by integrating multi-carrier energy devices and different energy distribution infrastructures under the framework of the smart micro energy internet, the energy management of such systems is challenging due to multiple decision makers having competitive or/and cooperative targets and operating in stochastic conditions and with asymmetric market information [21]. In these scenarios, traditional single-agent based control and decision theory is inadequate. As such, an advanced modeling, analysis, and decision tool is necessary for energy management in a smart micro energy internet.

Game theory has been a fundamental mathematical tool of economists, politicians, and sociologists for decades due to its capability of solving decision-making problems involving multiple agents with cooperative or/and competitive goals [22]. The concepts, theories, and methodologies of game theory can be used to guide the resolution of engineering design, operation, and control problems in a canonical and systematic way [23], [24]. Game theory has been widely used in smart grid planning, economic dispatch, and market equilibrium analysis [24]–[27]. Specifically, engineering game theory has become a powerful tool for the development of advanced dispatch and control schemes for the smart micro energy internet. Thus, in this paper, engineering game-theory based dispatch and control methods are investigated to help in realizing a self-approaching-optimum energy management system of the smart micro energy internet.

When developing advanced energy infrastructures, the construction of engineering demonstration sites is important for testing and validation. To this end, a prototype system of a solar-based smart micro energy internet has been established under the guidance of the smart micro energy internet framework in Qinghai University, China. The system is composed as follows: three kinds of solar-based power and heat sources, including a multi-function PV station, a solar chimney power station, and a full solar spectrum power generation; an energy hub based on solar-thermal compressed air energy storage system (ST-CAES); a smart micro-grid for library with building integrated photovoltaics; electricity demands such as carbon fiber recycling system and campus load; multi-carrier demands (heating, cooling, and power), such as solar-based wooden house; and an engineering game theory based self-approaching-optimum energy management system.

The rest of this paper is organized as follows. Section II proposes the vision of a smart micro energy internet in a comparative framework, while establishing concepts and design principles of the smart micro energy internet. A multi-criteria self-approaching-optimum energy management scheme based on engineering game theory for a smart micro energy internet is developed in Section III. A conceptual solar-based smart micro energy internet prototype system is presented in Section IV. Conclusions are drawn in Section V, and the intriguing research directions are also envisioned.

II. TRANSITION TO THE SMART MICRO ENERGY INTERNET

Solutions such as smart grid, integrated energy system, and energy internet have all been proposed in recent years to satisfy increasing energy demands in a secure and reliable way. However, most of these concepts focus on the transmission level side with considerably less attention given to the distribution level and demand side of energy systems. Taking into consideration the high-efficiency and flexibility characteristics of distributed co-generation, we focus here on the distribution and demand side as a way to accommodate micro-grids, distributed poly-generation systems, energy storage facilities, and associated energy distribution infrastructures within the framework of the smart micro energy internet. The next section is dedicated to exploiting and articulating concepts, design principles, and a vision of the smart micro energy internet.

*A. From Smart Grid to Smart Micro-grid*

As a typical next generation electric power system, the smart grid has been studied for decades. The smart grid is a modern electric power grid infrastructure that is efficient, reliable, and safe, offering seamless integration of renewable and alternative energy sources through automated control and modern communications technologies [28], [29]. To simplify the implementation of smart grid functions, such as reliability, self-healing, and load control, the concept of micro-grid has been proposed [19].

The micro-grid is regarded as the distributed way of the smart grid, having two fundamental operation modes: the autonomous mode and the grid-connection mode [30]. The micro-grid pays more attention to the distribution level and demand side. The generated power in a micro-grid is mainly for self-use in the autonomous mode. Insufficient or surplus power can be regulated with the connected PDN or power utility in the grid-connection operation mode.

The micro-grid with self-approaching-optimum energy management and satisfied multi-criteria performance of safety, economy, and environment protection can be treated as a smart micro-grid. A self-optimum-approaching dispatch and control center is necessary for the energy management system of the smart micro-grid. Advanced energy management strategies for both smart grid and micro-grid have been investigated to exploit their functionalities [31]-[33]. With the popularity of distributed CHP and micro-CCHP, the smart micro-grid is now capable of supplying cooling, heating, and power for a large portion of residential, industrial, and commercial demands, thereby yielding a CCHP micro-grid [34]. Similar to a traditional micro-grid, the CCHP micro-grid also requires flexible operations in both the autonomous mode and the energy distribution infrastructure connected mode. To facilitate the plug and play of such CCHP micro-grids requires a smart distribution-level energy infrastructure. The smart micro energy internet, therefore, is a perfect framework to implement this function.

*B. From Energy Internet to Smart Micro Energy Internet*

The energy internet is an energy management system with the electrical system as its core and the smart grid as its basis. It incorporates renewables and accommodates multi-carrier energy capabilities for improving the overall energy utilization ratio with advanced information, communications technologies, and power electronic technology [35], [36]. Inside the energy internet is the energy hub, which is a key facility wherein multiple energy carriers can be converted, conditioned, and stored. Specifically, the energy hub represents an interface between different energy infrastructures and/or load in the energy internet. Energy hubs consume power at their input ports, e.g., electricity and natural gas infrastructures, and provide certain required energy services such as electricity, heating, and cooling at the output ports. Although the energy hub connects multiple infrastructures in an energy internet, the complexity of the energy internet puts barriers on achieving some of its predefined functions.

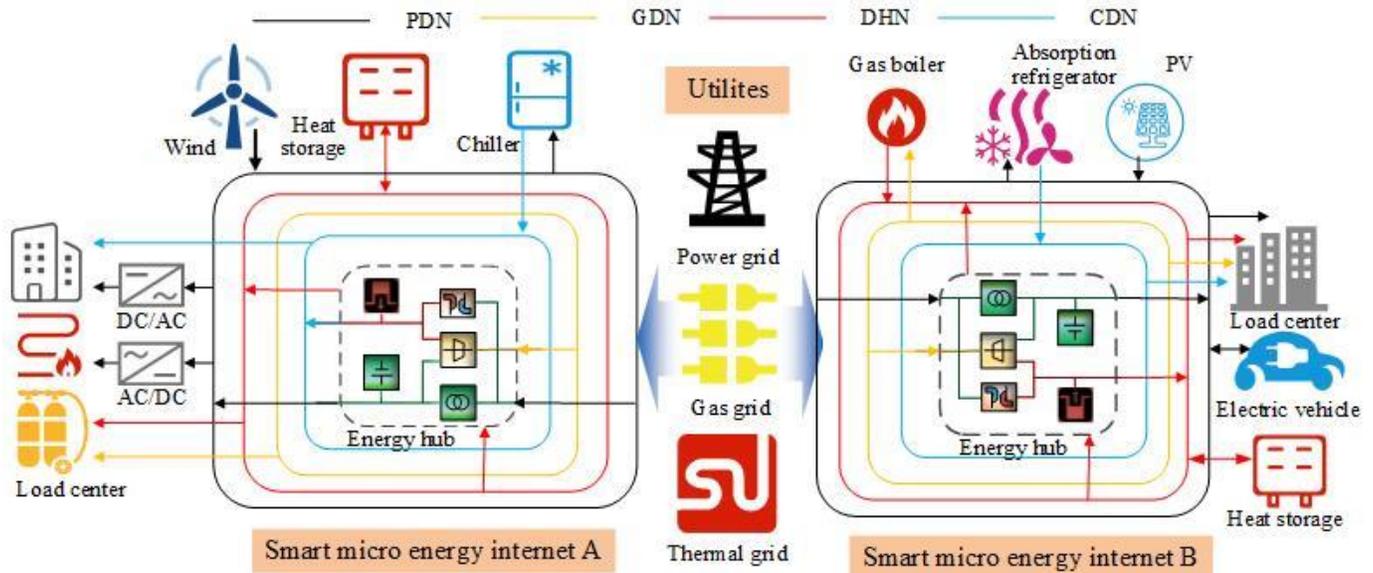

Fig. 1. Configuration of the networked smart micro energy internets.

To facilitate the functions of the energy internet at the distribution level, we investigate the smart micro energy internet as a class of energy networks composed of PDN, DHN, and CDN with energy hubs. The architecture of networked smart micro energy internets is depicted in Fig. 1. Energy infrastructures including electricity, district heating, natural gas, and transportation network are coupled through energy hubs, electric vehicle, heat storage system, and absorption refrigerator. Fig. 1 shows that the networked smart micro energy internets are connected to the energy utilities of natural gas, power, and thermal grid. Inside a smart micro energy internet, several energy distribution infrastructures are integrated via energy hubs, distributed poly-generation systems, or energy storage systems. Specifically, in the smart micro energy internet B, the interconnections among different distribution infrastructures can be clarified as: The energy hub builds the linkage among GDN, PDN, DHN, CDN；2）The absorption refrigerator connects PDN and CDN; 3) The gas boiler associates GDN and DHN. To design a smart micro energy internet, several fundamental principles are necessary, and are as follows:

1) At least one clean energy resource is integrated.
2) Energy storage facility with sufficient capacity is deployed.
3) Electricity can be distributed with the form of other energy carriers such as heating and cooling by using electric-heating and electric-cooling facilities.
4) Energy is mainly for self-use especially on demand side.
5) A multi-criteria self-approaching-optimum energy management system is incorporated to capture the unattended capability.

It should be mentioned that the "multi-criteria self-approaching-optimum" is a framework; it is composed of two fundamental aspects, i.e., 1) It is multi-criteria, which requires the designed system to be operated as one based on its operation targets. The definition of multi-criteria can be case by case, but as for the micro energy internet, it will include economic, security, and environmental considerations. 2) The second feature is that the framework is "self-approaching-optimum," which describes the ability of reaching the optimum operation point by itself, i.e., it is intelligent or smart. Definitely, the realization of "multi-criteria self-approaching-optimum" depends on mature optimization approaches, including both traditional methods as well as advanced methods. In this respect, the multi-criteria self-approaching-optimum is a good choice for the operation of the smart micro energy internet.

As indicated in Section II, a micro-grid is generally connected to the PDN, which in turn is integrated with the smart micro energy internet. Therefore, the concept of the micro-grid is within the scope of the smart micro energy internet. Under the vision of the energy internet, the grid-connection mode of the micro-grid means that it is connected to the smart micro energy internet. Since the smart micro energy internet can accommodate the traditional smart micro-grid, the voltage level and electricity capacity/demand is similar to that of the smart micro grid. As for the heating network, since most of the existing district heating system's piping network is distributed locally, it is natural to be integrated to the proposed smart micro energy internet, especially for the 4th generation of district heating and cooling techniques [13]. The key parameters in the district heating system are the temperature at the supply and return systems, as well as the hydraulic distribution along the piping networks. The interconnection points between the electric part and heat part of the smart micro energy internet are the co-generation units as well as the energy hubs. The smart micro energy internet is integrated with the energy utility infrastructure, including power, gas, and thermal grids. The local gas distribution network can also be incorporated in the smart micro energy internet as shown in Fig. 1. In this case, the interconnections among distribution level infrastructures will be realized via the gas-fired CHP, gas-fired CCHP, in addition to the energy hubs.

We firmly believe that the implementation of the smart micro energy internet can simplify energy internet functionalities, such as multi-carrier energy synergy that incorporates more renewable energy resources with improved reliability. Undoubtedly, the performance of the smart micro energy internet depends on its dispatch and control strategies. A multi-criteria, self-approaching-optimum operation energy management system is required to realize the secure, cost-effective, and environmentally friendly operation of the smart micro energy internet [20].

## III. ENGINEERING GAME-THEORY METHODS FOR THE SMART MICRO ENERGY INTERNET

Smart micro energy internet is typically characterized by incorporating multiple energy sources (wind, PV, natural gas suppliers, heat suppliers, etc.), multi-network (power network, gas network, heating network, cooling network), multi-user (power users, heat users, natural gas users, central air conditioning systems, etc.) and multiple energy forms (electricity, gas, heat, cooling, etc.). Thus, the smart micro energy internet is a complex stochastic system integrated laterally by electricity, gas, heat, cold, and other energy distribution infrastructures, and composed longitudinally by different decision-making agents including energy suppliers, dispatch centers, and multiple users. In this regard, the planning and energy management of the smart micro energy internet calls for advanced methods.

As indicated previously, the emphasis of game theory is on the decision behaviors of individual players with competing or cooperating goals based on the collected information. The focus of engineering game theory is on the application of fundamental concepts, modeling, and solution methods of game theory to the decision problem of engineering design and experiment; it takes into account engineering conditions [24]. Motivated by [40], we develop here engineering game-theory methods for the smart micro energy internet as illustrated next.

### A. General Remarks

The engineering game theory approach originated from our previous work in the field of smart grid [23], [24]. The fundamental motivation was to convert the decision and control problems in engineering to mathematical problems using game modeling techniques. These problems could then be solved using equilibrium methods that provide the engineering decision-maker a competitive and optimized solution. Different from traditional game-theory methods, we put more emphasis on the competitive status among decision players, whereas engineering game-theory methods pay more attention to the competitive behavior among humans and the nature. In this manner, we bridge traditional game-theory methods to engineering applications, and in this instance specifically to planning and energy management of the smart micro energy internet. The fundamental motivations of the engineering game-theory methods can be illustrated as:

$$\min_{x_1 \in X_1} f_1(x_1, x_2) \min_{x_2 \in X_2} f_2(x_1, x_2) \tag{1}$$

where nature $x_1$ is a fictitious player, with rationality $f_1(x_1, x_2)$, and strategy set $X_1$ while $x_2$ denotes the decision maker with the utility function $f_2(x_1, x_2)$ and strategy set $X_2$. Analyzing nature's rationality $x_1$ and modeling the strategy set $X_1$ are both difficult propositions. The basic framework of the engineering game-theory method is depicted in Table I.

TABLE I
GENERAL REMARKS ON ENGINEERING GAME-THEORY METHODS

| Decision Problem | Game Model | Solution Method |
|---|---|---|
| Multi-objective optimization | Stationary game | Nash equilibrium |
| Robust control | Differential game | Feedback Nash equilibrium |
| Robust optimization | Zero-sum game | Saddle Nash equilibrium |
| Multi-level optimization | Leader-follower game | Nash-Stackelberg-Nash equilibrium |
| Security defense | Security game | Bayes-Nash equilibrium |
| System evolution | Evolution game | Evolutionary stability equilibrium |

### B. Planning of the Smart Micro Energy Internet

The planning of the smart micro energy internet mainly focuses on the optimal allocation and size of energy hubs, distributed generators, etc. Since power distribution, natural gas transportation, and district heating infrastructures are existing entities, planning of the smart micro energy internet here is focused on the energy hubs. Specifically, the energy hubs in a smart micro energy internet are required to capture the features of combined heat, power, cooling generation, and storage. Thus, an energy hub is usually an integration of different components, including transformer, micro-turbine, heat exchanger, furnace, absorption chiller, power energy storage unit, and hot water storage.

For the sake of feasibility, an optimal selection of these components is essential for the smart micro energy internet as in [43]. Cost saving and emissions reduction are the fundamental criteria for the smart micro energy internet. In this perspective, the optimal allocation and size of the energy hub is multi-criteria, and thus multi-objective programming methods are usually adopted to handle this kind of problem. In general, it is difficult to reach the two minimums simultaneously due to the contradiction between carbon emission and system cost. Therefore, there is a need to develop evolutionary algorithms (EA) that is capable of providing many non-inferior solutions to tackle this problem. EA methods have been viewed as effective in obtaining the Pareto front. However, in an actual smart micro-energy internet, the decision-maker needs one balanced solution, and as such, it will be difficult to select one from the Pareto front.

Stationary game based multi-objective optimization methods can be adopted to determine a balanced solution from the Pareto front. For this purpose, three methods have been proposed: the comprehensive method, the weighting method, and the Nash bargaining method [24]. Readers are referred to this work for additional details. In this subsection, we focus on the Nash bargaining method of the most commonly bi-objective optimization. Specifically, the planning problem in terms of cost saving and emission reduction can be generalized as:

$$\min\{f_1(x), f_2(x)\}$$
$$\text{s.t. } g(x) = 0, h(x) \leq 0, x \in R^m \quad (2)$$

where $f_1(x)$, $f_2(x)$, respectively, denote the cost and carbon emission while $g(x)$ and $h(x)$ are the constraints.

From the perspective of engineering game theory, the two objectives in (2) are treated as two virtual players who negotiate with each other on how to distribute the planning resources. Based on Nash's theory, the bargaining solution of (2) is the optimal solution of (3):

$$\min_{x_1 \in X_P} \left(f_1^d - f_1(x)\right)\left(f_2^d - f_2(x)\right) \quad (3)$$

where $f_1^d$ and $f_2^d$ are the maximum cost and carbon emission. The problem (3) can be solved through a univariate parametric method such as the golden section search algorithm. The Nash bargaining method herein can be directly extended to the multi-objective one. Another advantage of the Nash bargaining method is that a no weight parameter is needed to normalize the magnitudes of the two objectives, i.e., it remains invariant under linear transformation of the objective function.

### C. Hierarchical Energy Management of the Smart Micro Energy Internet

The coupling among multiple physical systems, multiple decision-making agents, and the stochastic market environment makes the dispatch and control of the smart micro energy internet a challenging issue. In this regard, to achieve multi-criteria self-approaching-optimum dispatch and control for the smart micro energy internet, we propose an engineering game-theory based hierarchical energy management scheme in this subsection. The comprehensive energy management system consists of three layers: utility connection layer, intra energy net layer, and component layer as illustrated in Fig. 2.

#### 1) Utility Connection Layer

This layer decides the operation mode of smart micro energy internet according to the prices of electricity and fuels in the energy market and other information, such as power demand and weather prediction. For the grid-connection mode, the exchanged multi-carrier flow among networked smart micro energy internets and the connected energy utilities, as well as the exchanged multi-carrier energy with other smart micro energy internets should be determined. In the case of the autonomous mode, the exchanged multi-carrier flows among micro energy internets, and energy utilities are set to zero. Since there are multiple decision makers in this layer, the difficulty is in determining the exchanged multi-carrier flow among micro energy internets due to inevitable competitive or/and cooperative targets of decision makers, who operate under the stochastic conditions and asymmetric market information [38], [39].

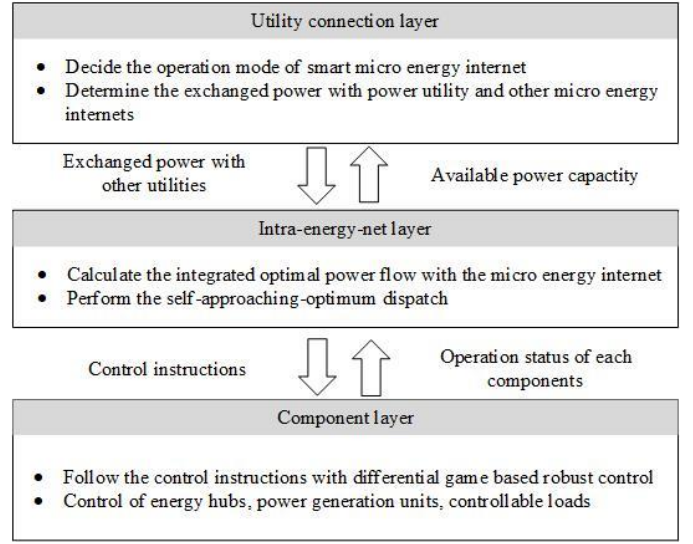

Fig. 1 Hierarchical energy management of smart micro energy internet

To this end, engineering game-theory based methods can be adopted while at the same time taking into account the utility of different market participants in terms of safety, economy, and environmental protection criteria [40]. The decision behavior among the smart micro energy internets can be formulated as a Nash game in terms of cooperative or non-cooperative targets. Then, the exchanged power can be determined by a Nash equilibrium, as suggested in Fig. 3. Once the exchanged power flow is regulated, the integrated optimal power flow (IOPF) among multi-carrier energy utilities within the smart micro energy internet can be determined by the intra-energy-net layer.

#### 2) Intra-energy-net Layer

Optimal power flow among multi-carriers is needed for integrated energy infrastructures [41]. The intra-energy-net layer is responsible for the determination of IOPF within the smart micro energy internet once the exchanged power is calculated in the utility connection layer. According to the features of multi-carrier energy infrastructures, the IOPF of the smart micro energy internet usually contains multi-decision makers and has multi-time scale characteristics. The determination of IOPF is also daunting since the multi-decision makers may have a competitive or cooperative utility function.

When it comes to a specific smart micro energy internet incorporating natural gas, district heating, and power distribution network with CHP, three types of dynamics exist at different time-scales: thermal dynamic, gas, and electrical dynamic. To overcome these multi-time-scale characteristics, the intra-energy-net layer can be further divided into three sub-layers to handle slow, medium, and fast control layer as in [42]. To address the feature of multiple decision makers, which is important since infrastructures are usually owned by different utilities in today's energy market, two situations need to be considered. In the first scenario, different infrastructures have the same decision maker. This can be solved through a cooperative game. In the second scenario, different infrastructures have disparate decision makers. This case can be handled by the non-cooperative game, as depicted in Fig. 3, thus yielding a leader-follower game.

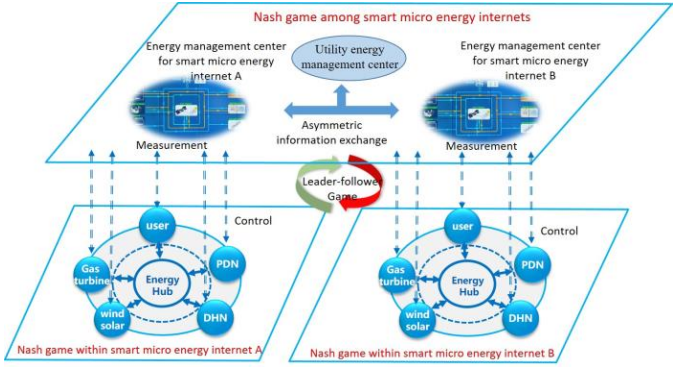

Fig. 2 Engineering game-theory methods based control and dispatch strategy for smart micro energy internet.

Once the IOPF within the smart micro energy internet is determined, dispatch strategies and instructions can be executed through the elements in each component layer, including the distributed poly-generation unit, the energy hubs, and the energy storage facilities.

*3) Component Layer*

In the smart micro energy internet, the component layer is responsible for real-time control of the device for executing the setting points from the intra-energy-net layer. These devices usually suffer from external disturbance such as measurement noise, and internal disturbance such as unmodeled dynamics. Therefore, the control strategies for the devices in this layer need to be robust to follow the setting values.

The differential game based robust control is a proper candidate for this scenario; it treats uncertainty as a fictitious player competing with a human player to decide the performance of the devices. Mathematically, the controlled device in the smart micro energy internet can be modeled as:

$$\begin{cases} \dot{x} = f(x) + g_1(x)w + g_2(x)u \\ z = h(x) + k(x)u \end{cases} \quad (4)$$

where $f_1(x), g_1(x), g_2(x), h(x), k(x)$ are the device dynamics, $x$ is the state variable, $u$ the control variable, and $w$ models the uncertainty. The human decision-maker has the goal of minimizing the utility function:

$$J(u,w) = \int_0^T \left( \|z\|^2 - \gamma^2 \|w\|^2 \right) dt \le 0, \forall T \ge 0 \quad (5)$$

While the uncertainty wants to maximize $J(u, w)$, we get the differential game model of robust control:

$$\min_u \max_w J(u,w) \le 0$$
$$s.t. \quad \dot{x} = f(x) + g_1(x)w + g_2(x)u \quad (6)$$

The control law can be obtained by solving the differential game and as a means for finding the feedback Nash equilibrium. Several kinds of feedback Nash equilibrium seeking methods, including variable metric feedback linearization, controlled Hamilton designing method, policy iteration method, and approximated dynamic programming method are summarized in [24].

## IV. ENGINEERING IMPLEMENTATION: SOLAR BASED SMART MICRO ENERGY INTERNET

In addition to developing an energy management scheme, the construction of the smart micro energy internet itself is also critical. To this end, the first domestic implementation of the conceptual solar-based smart micro energy internet has been undertaken in Qinghai University. The motivation was to develop a totally clean smart micro energy internet with advanced solar energy utilization techniques. The fundamental idea was to generate electricity, thermal energy, and cold energy with a series of PV and solar thermal techniques to satisfy the multi-carrier demands (heating, cooling, and power) in the smart micro energy internet. The generated electricity can be either for self-use or being fed into the power utility to which the smart micro energy internet is connected. Heat collected by the parabolic trough collectors can be either used to generate electricity through a turbine within the compressed air energy storage system or to be directly fed into the DHN. Cold energy is produced by absorption refrigerator using exhaust heat in the system and injected into the CDN. A solar based wooden house is used to simulate the multi-carrier energy demands in the smart micro energy internet. Generally, the power, heating, and cooling demands can be fulfilled by the smart micro energy internet itself. For specific cases, the electricity demands can be supplied by the power utility to which the solar-based smart micro energy internet is connected.

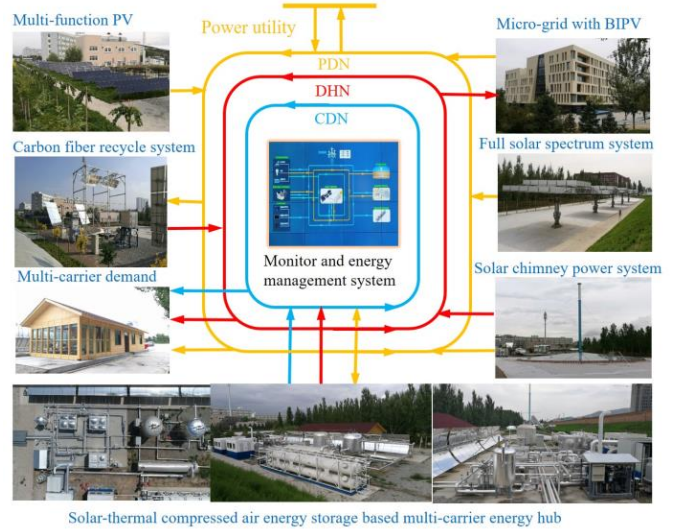

Fig. 4. Layout and multi-carrier energy flow of the solar-based smart micro energy internet in Qinghai University.

The constructed solar-based smart micro energy internet is composed of three kinds of solar-based power and heat sources: multi-function PV station, solar chimney power station, and full solar spectrum power generation. It also contains an energy hub in the form of a solar-thermal compressed air energy storage facility, a smart micro-grid for the library building integrated with photovoltaics (BIPV), multi-carrier loads including carbon fiber recycling system, solar-based wooden house, and campus load, and the engineering game-theory methods based multi-criteria self-approaching-optimum energy management system. The layout and the multi-carrier energy flow of the smart micro energy internet are shown in Fig. 4. Note that, although the construction of this smart micro energy internet was completed in August 2016, it was still evolving to integrate more and more infrastructures to target the internet of things. Each subsystem and the energy management center of the solar-based smart micro energy internet are elaborated as follows.

*A. Energy Hub: ST-CAES*

The solar-thermal compressed air energy storage system (CAES) is a typical CAES facility, which captures the features of zero-carbon emission, collecting and recycling of air compression heat and solar thermal, and combined heat and power energy storage, and tri-generation. It creates the link between the PDN, DHN, as well as the CDN in the constructed solar-based smart micro energy internet. The ST-CAES acts as an energy hub and its architecture is shown in Fig. 5.

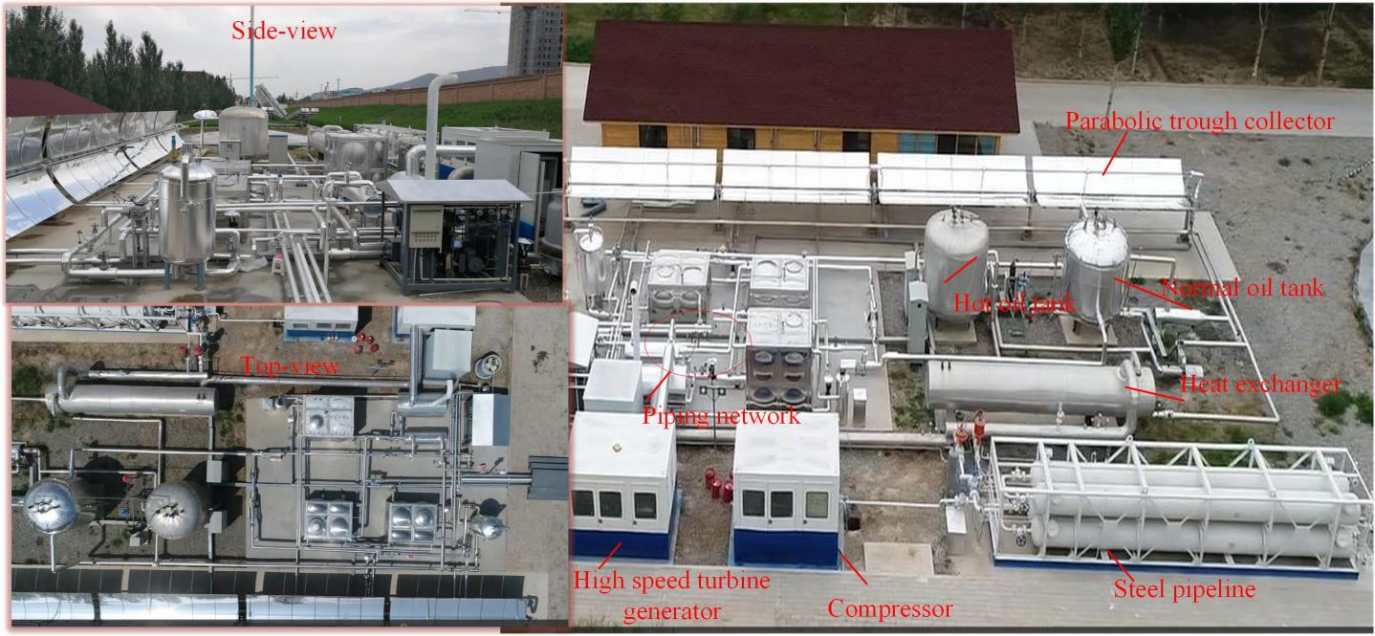

Fig. 5. Architecture of the solar-thermal compressed air energy storage system based energy hub.

The operation of the ST-CAES facility is comprised of two fundamental processes. During the charging process, electricity, generated by the power sources in the smart micro energy internet, is utilized to drive the compressor to generate high-pressure air stored in a steel pipeline. Meanwhile, the solar thermal energy collected using parabolic trough collectors (PTCs) during good weather, as well as the heat generated, along with the compressed air are stored in the thermal storage tank. During the discharging process, the stored high-pressure air is preheated with the stored thermal energy to drive the high-speed turbine generator to produce electricity. The utilization of collected solar thermal energy improves the temperature of the suction air of the high-speed turbine generator and also increases the efficiency of the ST-CAES. The surplus thermal energy can be injected into the DHN. Moreover, the exhaust air of the turbine can be used for cooling by regulating the temperature of the suction air. The exhaust heat is utilized to generate cooling energy if needed through the absorption refrigerator equipped with the ST-CAES. The generated cooling energy can be injected into CDN. Thus, the ST-CAES allows the transfer of power, heat, and cooling in space and time and is typically in the form of a combined cooling, heating, and power energy storage system. The detailed dispatch model of the CAES hub and its effectiveness in a typical micro energy internet contains PDN and DHN are investigated in [16].

*B. Solar-thermal Based Carbon Fiber Recycling System*

In the constructed solar-based smart micro energy internet, a solar-thermal based carbon fiber recycling system is deployed as the power load. The architecture of this system is depicted in Fig. 6. The system is composed of controllable solar reflectors, concentrators, pre-treatment platform, fine treatment furnace, and solar thermal storage tank. The electricity is needed to drive the solar reflectors, concentrators, etc. It is well to be reminded that the fine furnace uses molten salt as the heat transfer material whose temperature can reach 600 °C with a ±2°C tolerance; thus technically this system can be regarded as a class of solar thermal power station. In short, although the system is used for recycling carbon fiber material, its function is equivalent to a solar thermal power plant. In this regard, this subsystem can be viewed as either a power load or generator with combined heat and power generation capability, and its model is similar to the CSP in [17].

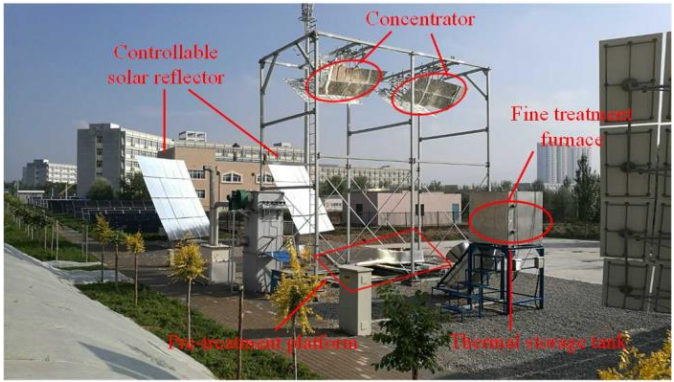

Fig. 6. Solar-thermal based carbon fiber recycling system.

*C. Power Sources in the Smart Micro Energy Internet*

Since the essence of the solar-based smart micro energy internet is to utilize solar energy as the source for each energy carrier, power sources are all based on existing solar utilization techniques, including photoelectric conversion, thermoelectric conversion, and combined photoelectric and thermoelectric techniques. The representatives of such techniques in the constructed smart micro energy internet are respectively the multi-function PV station, the plateau solar chimney power station, and the full solar spectrum power station, as shown in Fig. 7 through Fig. 9.

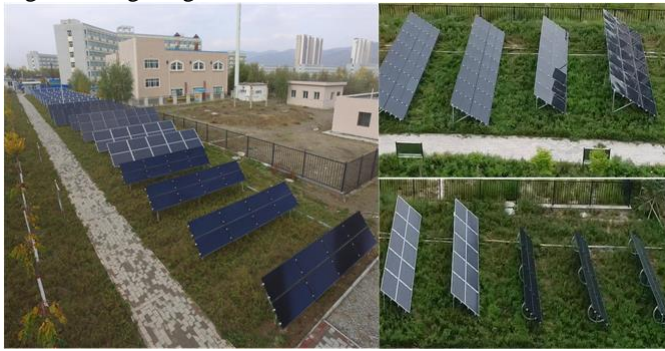

Fig. 7. Multi-function PV station.

*1) Multi-function PV Station*

The multi-function PV station is one of the main power sources for the smart micro energy internet. Its annual power generation capacity can reach 80,000 kWh. The generated electric power is mainly consumed by the university campus. The PV station consists of 6 kinds of photovoltaic cell assemblies, including monocrystalline, polycrystalline, double-sided, cadmium telluride, amorphous, and copper indium gallium selenide, three types of inverters, micro-inverters, intelligent optimizer, and string inverter; and 10 combinations incorporating different modules, inverters, and brackets.

*2) Plateau Solar Chimney Power Station*

Representing the thermoelectric technique, the solar chimney power generation system is mainly composed of heat-collecting shed, turbine, chimney, and heat storage layer. Solar thermal is utilized with the purpose of enhancing the air temperature of the heat-collecting shed, reducing air density, and generating pressure differences with the ambient environment. The air is then drawn by the chimney to generate an updraft air flow to drive the turbine to produce electricity. This subsystem features no pollution during the operation, can be easily set up in the desert; it provides local environment improvement by minimizing greenhouse effects.

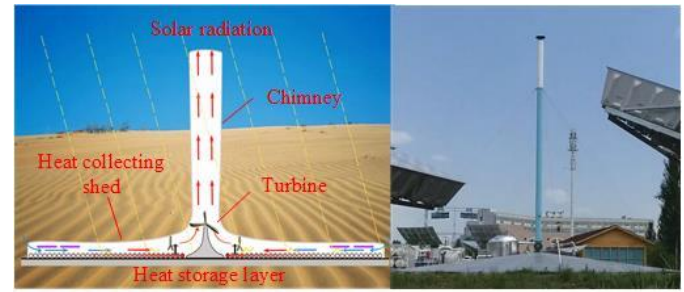

Fig. 8. Plateau solar chimney power generation.

*3) Full Solar Spectrum Power Generation*

Generally, 99% of the solar radiation energy is located in the wavelength range of 200 to 3000 nm. 58% of this solar energy is distributed in the form of ultraviolet and visible light with 200–800 nm wavelength, while 42% energy is in infrared light with the wavelength coverage of 800–3000 nm. In this regard, the photoelectric technique and thermoelectric technique individually utilize 58% and 42% of the solar radiation energy to produce electricity. However, since these two techniques rely heavily on the conversion material, the efficiency is usually low due to the fact that only the partial spectral range of the solar radiation energy is utilized. Full solar spectrum power generation technique can fix this issue. The essence of the full solar spectrum power generation system is to integrate solar thermoelectric conversion techniques and photoelectric conversion techniques to yield a high efficiency combined thermoelectric and photoelectric power generation technique. The fundamental principles of full solar spectrum power generation technique and a schematic of the system is shown in Fig. 9.

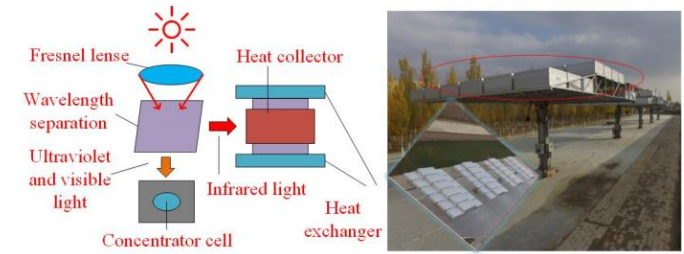

Fig. 9. Full solar spectrum power generation system.

*D. Solar-based Micro-grid with BIPV Technique*

Building integrated photovoltaics is seen as one of the main developmental trends in urban buildings and energy systems since they offer the most efficient way to implement distributed PV generation and micro-grid technologies. BIPV combines PV and the building, with no extra floor space needed for installment; this type of local accommodation is the best application form of PV in urban areas, especially, in the far-western regions of China. In this respect, a solar-based micro-grid with BIPV technique is built for the university library as shown in Fig. 10.

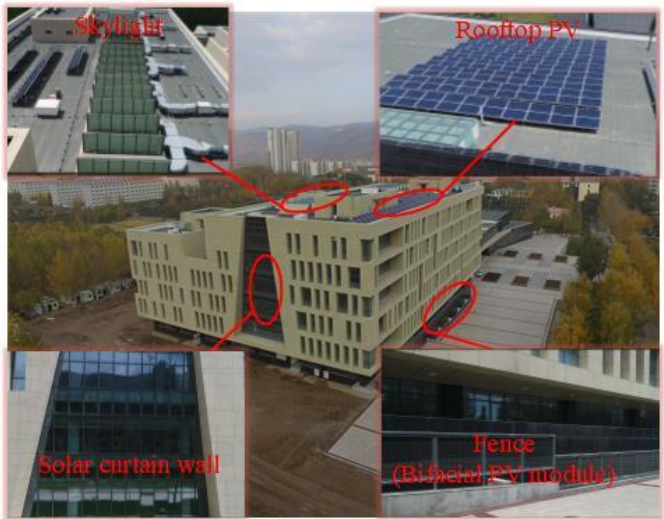

Fig. 10. Solar-based micro-grid with BIPV technique.

The solar-based micro-grid is the first BIPV innovation platform in Qinghai province. It has facilitated the coordination between green energy technologies and modern architectural aesthetics. This micro-grid is actually a CCHP micro-grid and it is connected to the solar based smart micro energy internet. Thus, the library micro-grid can act as an energy source as well as energy load, depending on its operation mode.

*E. Energy Management System*

Energy management system (EMS) of the solar-based smart micro energy internet has been developed to realize smart operations. The developed EMS uses the collected operation measurements including multi-carrier energy generation and demand, voltage level in the PDN, temperature and pressure of piping system in DHN and CDN, and real-time operation data to dispatch and control each facility in the smart micro energy internet in a centralized manner.

Since there is no other smart micro energy internet, the utility connection layer of EMS decides the operation mode of this micro energy internet, while calculating its unbalanced energy demand. The intra-energy-net layer guarantees optimal energy distribution among the generation sources and multi-carrier demands. The component layer follows the dispatch instructions from the intra-energy-net layer to robustly control the PV station operations, full solar spectrum system, solar chimney power station, energy hub, and the micro-grid. Up until April 2017, the accumulated electricity generated by multi-function PV station and BIPV based micro-grid had reached 109.91 MWh and 20.88 MWh, respectively.

## V. CONCLUSION

This paper proposes the concept and design principles of the smart micro energy internet to implement the vision of an energy internet. An engineering game-theory based energy management scheme is developed to realize the multi-criteria self-approaching-optimum operation of the smart micro energy internet. A prototype of a solar-based conceptual smart micro energy internet system has been established in Qinghai University. The micro energy internet is the distributed form of the energy internet. The use of the micro energy internet helps to simplify the implementation of many energy internet functions. The smart micro energy internet is capable of multi-criteria self-approaching-optimum operation. Engineering game-theory methods provide effective ways to realize the real-time self-approaching-optimum dispatch and control strategies.

The vision of the smart micro energy internet relies on substantial advances in intelligent distributed or decentralized control and decision mechanisms. Our future work will focus on the development of engineering game-theory based distributed control and optimization schemes as well as their implementation in operating the smart micro energy internet.

ACKNOWLEDGMENT

The authors would like to thank Dr. Shaowei Huang, Tsinghua University, Mr. Yang Si, Qinghai University, for their excellent work in the construction of the solar based smart micro energy internet, and Prof. Feng Liu, Prof. Wei Wei, Tsinghua University, for their valuable insights in the implementation of game-theoretic methods to the smart micro energy internet.